\documentclass{elsart3-1}

\usepackage{amsmath,amssymb}
\usepackage[english,francais]{babel}
\usepackage{xcolor}

\newtheorem{Theorem}{Theorem}[section]
\newtheorem{Lemma}[Theorem]{Lemma}

\newtheorem{Corollary}[Theorem]{Corollary}

\newtheorem{Theoreme}{Th\'eor\`eme}

\newtheorem{Corollaire}[Theorem]{Corollaire}

\renewcommand{\theequation}{\arabic{equation}}
\DeclareTextSymbol{\degre}{OT1}{23}

\def	\N	{	\ensuremath	{\mathbb	{N}	}	}
\def	\Q	{	\ensuremath	{\mathbb	{Q}	}	}
\def	\Z	{	\ensuremath	{\mathbb	{Z}	}	}

\def	\C	{	\ensuremath	{\mathbb	{C}	}	}

\newcommand {\seq}[1] {	\ensuremath { \underline{\mathbf{#1}\!}\, } }

\journal{the Acad\'emie des sciences}

\begin{document}

\centerline{}

\begin{frontmatter}

\selectlanguage{english}

\title{The Algebra of Hurwitz Multizeta Functions}
\author{Olivier Bouillot}
\ead{olivier.bouillot@univ-mlv.fr}
\ead[url]{http://www-igm.univ-mlv.fr/$\sim$bouillot/}
\address	{	Laboratoire d'informatique Gaspard-Monge (LIGM)	\\
					Universit\'e Paris-Est Marne-la-Vall\'ee			\\
					Cit\'e Descartes - $5$ bd Descartes - 77454 Marne-la-Vallée Cedex 2		\\
				}

\medskip

\begin{center}
    \small Received *****;
     accepted after revision +++++\\
    Presented by £££££
\end{center}

\begin{abstract}
    \selectlanguage{english}
    Multizeta values are real numbers which span a complicated algebra:
    there exist two different interacting products. A functional analog of these numbers is defined so as
    to obtain a better understanding of them, the Hurwitz multizeta functions, which span an algebra for
    which a precise description is wanted.
    In this note, we prove that the algebra of Hurwitz multizeta functions is a polynomial algebra.
    

    \vskip 0.5\baselineskip

    \selectlanguage{francais}
    \noindent{\bf R\'esum\'e} \vskip 0.5\baselineskip \noindent
    {\bf Sur l'alg\`ebre des multiz\^etas de Hurwitz. }
    
    Les multiz\^etas sont des nombres r\'eels poss\'edant une structure d'alg\`ebre
    complexe : il existe deux produits interagissant.
    Il est naturel de d\'efinir un analogue fonctionnel de ces nombres pour en avoir une meilleure
    compr\'ehension, ce qui conduit aux multiz\^etas de Hurwitz, dont on souhaiterait
    conna\^itre pr\'ecis\'emment la structure d'alg\`ebre.
    Dans cette note, nous montrons que l'alg\`ebre des multiz\^etas de Hurwitz est une alg\`ebre
    de polyn\^omes.
\end{abstract}

\end{frontmatter}


\selectlanguage{french}
\vspace{-1cm}
\section*{Version fran\c{c}aise abr\'eg\'ee}

Pour toute s\'equence $(s_1, \cdots, s_r) \in (\N^*)^r$, de longueur $r \in \N^*$ telle que $s_1 \geq 2$,
on appelle multiz\^eta de Hurwitz (ou encore polyz\^eta monocentr\'e de Hurwitz) la fonction d\'efinie sur
$\C - (-\N)$, par :
\begin{equation}
	\mathcal{H}e^{s_1, \cdots, s_r} : z \longmapsto \sum	_{0 < n_r < \cdots < n_1}
								\frac	{1}
									{	(n_1 + z)^{s_1}
										\cdots
										(n_r + z)^{s_r}
									}	\ .
\end{equation}

La condition $s_1 \geq 2$ assure la convergence de la s\'erie ainsi que le caract\`ere holomorphe de
$\mathcal{H}e^{s_1, \cdots, s_r}$ sur $\C - (-\N)$~. On d\'efinit aussi 
$\mathcal{H}e^\emptyset : z \longmapsto 1$ (Cf. \cite{Bouillot-Thesis} et \cite{Bouillot-MTGF} pour un
contexte o\`u ces fonctions apparaissent naturellement en dynamique holomorphe, en lien avec les
multitangentes ; Cf. \cite{Minh}, \cite{Minh4} et \cite{Murty-Sinha} pour d'autres travaux se rapportant \`a ces
fonctions)~.

Une telle fonction est une \'evaluation d'une fonction quasi-sym\'etrique monomiale
$M_{s_1, \cdots, s_r}$~(Cf. \cite{Gessel} pour la premi\`ere pr\'esentation de ces fonctions et
\cite{Bergeron} pour une pr\'esentation plus r\'ecente) . Ainsi, cette famille de fonctions engendre une
alg\`ebre not\'ee
$\mathcal{H}mzf_{cv}$, son produit \'etant le produit de $\textit{QSym}$, appel\'e aussi bien
stuffle, shuffle augment\'e, produit de battage contractant, $\cdots$, suivant les auteurs
(cf. \cite{Hoffman}, \cite{Ebrahimi-FardGuo}, \cite{Bouillot-MTGF}). L'objectif de cette note est de 
d\'ecrire la structure alg\'ebrique de $\mathcal{H}mzf_{cv}$~.

\bigskip

On peut remarquer que la famille des multiz\^etas de Hurwitz est un analogue fonctionnel de celle des
multiz\^etas,ces derniers \'etant d\'efinis par
$\mathcal{Z}e^{s_1, \cdots, s_r} = \mathcal{H}e^{s_1, \cdots, s_r}(0)$
(Cf. \cite{Cartier}, \cite{Waldschmidt} et \cite{Zudilin} pour une introduction aux multiz\^etas).
On sait peu de choses sur la nature arithm\'etique des multiz\^etas, malgr\'es quelques avanc\' ees
r\'ecentes (Cf. \cite{Brown1}, \cite{Brown2} , \cite{Minh2} et \cite{Minh3}).
Ces nombres engendrent aussi une alg\`ebre not\'ee $\mathcal{M}zv_{cv}$ et l'on peut demander une
description alg\'ebrique de celle-ci.
Par exemple, $\mathcal{M}zv_{cv}$ est-elle une alg\`ebre gradu\'ee, et si oui, quelle est la suite des
dimensions de ses composantes homog\`enes ?

D\'eterminer la structure alg\'ebrique d'une alg\`ebre non formelle, comme $\mathcal{M}zv_{cv}$ ou
$\mathcal{H}mzf_{cv}$, est en g\'en\'eral une question difficile car, \`a priori, d'autres produits
peuvent entrer en jeu et complexifier la structure de l'alg\`ebre \'etudi\'ee. C'est en particulier le cas
de $\mathcal{M}zv_{cv}$. Dans le cas des multiz\^etas de Hurwitz, la situation est en fait incomparablement
plus simple, car le stuffle, provenant de l'identification des multiz\^etas de Hurwitz avec une sp\'ecialisation 
des fonctions quasi-sym\'etriques monomiales, est le seul produit possible.

Dans cette note, nous montrons que $\mathcal{H}mzf_{cv}$ est en fait une alg\`ebre polynomiale. Le point
cl\'e est le%

\bigskip

\begin{Theoreme}
	La famille	$	\left (\mathcal{H}e^{s_1, \cdots, s_r} \right )
				_{	r \in \N
					\atop
					{	s_1, \cdots, s_r \in \N^*
						\atop
						s_1 \geq 2
					}
				}
			$ est $\C(z)$-lin\'eairement ind\'ependante.
\end{Theoreme}

\bigskip

Notons que la r\'ef\'erence \cite{Joyner} propose une preuve d'une propri\'et\'e plus faible, l'ind\'ependance lin\'eaire sur $\C$ des multiz\^etas de Hurwitz\footnote{Selon le rapporteur, qui nous a signal\'e cette r\'ef\'erence, cette preuve contient quelques erreurs...}.
Les corollaires suivants de ce th\'eor\`eme permettent de r\'epondre \`a toutes les questions que l'on peut poser concernant les multiz\^etas d'Hurwitz. Par exemple :

\bigskip

\begin{Corollaire}
	Notons $\textit{QSym}_{cv}$ la sous-alg\`ebre de $\textit{QSym}$ engendr\'ee par les mon\^omes
	$M_I$ o\`u $I = (i_1, \cdots, i_r)$ est une composition telle que $i_1 \geq 2$~.
	Alors :
	\begin{equation}
		\mathcal{H}mzf_{cv} \simeq \textit{QSym}_{cv}
				    \simeq
				    \Q	\big \langle
							\mathcal{L}yn(y_1 ; y_2 ; \cdots) - \{y_1\}
					\big \rangle
		\ ,
	\end{equation}
	o\`u $\mathcal{L}yn(y_1 ; y_2 ; \cdots)$ repr\'esente l'ensemble des mots de Lyndon sur l'alphabet
	$Y = \{ y_1 ; y_2 ; \cdots \}$.
\end{Corollaire}

Pour plus d'informations sur les mots de Lyndon, on renvoie le lecteur \`a \cite{Reutenauer}.

\bigskip

\begin{Corollaire}
	Toute relation alg\'ebrique de $\mathcal{H}mzv_{cv}$ est cons\'equence de l'\'evaluation du
	produit de stuffle.
\end{Corollaire}

\bigskip

\begin{Corollaire}
	Notons $\mathcal{H}mzv_{cv,n}$ l'alg\`ebre des multiz\^etas de Hurwitz de poids $n$, c'est-\`a-dire
	la sous-alg\`ebre de $\mathcal{H}mzv_{cv}$ engendr\'ee par les multiz\^etas de Hurwitz
	$\mathcal{H}e^{s_1, \cdots, s_r}$ v\'erifiant $s_1 + \cdots + s_r = n$.	\\
	Alors :
	\begin{enumerate}
	      \item $\mathcal{H}mzv_{cv}$ est gradu\'ee par le poids :
			\begin{equation}
				\left \{
					\begin{array}{l}
						\mathcal{H}mzv_{cv} = \bigoplus	_{n \in \N}
										\mathcal{H}mzv_{cv,n}
						\ .
						\\
						\forall (p , q) \in \N^2 \ , \ 
						\mathcal{H}mzv_{cv,p} \cdot \mathcal{H}mzv_{cv,q}
						\subset
						\mathcal{H}mzv_{cv,p + q} \ .
					\end{array}
				\right.
			\end{equation}
		\item $	\left \{
					\begin{array}{l}
						\text{dim } \mathcal{H}mzf_{cv,0} = 1
						\ , \ 
						\text{dim } \mathcal{H}mzf_{cv,0} = 0.
						\\
						\text{dim } \mathcal{H}mzf_{cv,n} = 2^{n - 1}
						\text{ for all }n \geq 2.
					\end{array}
				\right.
			$
 
	\end{enumerate}
\end{Corollaire}

\selectlanguage{english}

\section{Introduction}

For any sequence $(s_1, \cdots, s_r) \in (\N^*)^r$, of length $r \in \N^*$ such that $s_1 \geq 2$,
we define the Hurwitz multizeta function (or mono-center Hurwitz polyzeta), over $\C - (-\N)$, by:
\begin{equation}
	\mathcal{H}e^{s_1, \cdots, s_r} : z \longmapsto \sum	_{0 < n_r < \cdots < n_1}
								\frac	{1}
									{	(n_1 + z)^{s_1}
										\cdots
										(n_r + z)^{s_r}
									}	\ .
\end{equation}

The condition $s_1 \geq 2$ ensures the convergence and, consequently,
$\mathcal{H}e^{s_1, \cdots, s_r}$ is a holomorphic function over $\C - (-\N)$. We also define
$\mathcal{H}e^\emptyset : z \longmapsto 1$~. (See. \cite{Bouillot-Thesis} and \cite{Bouillot-MTGF} where these functions appear naturally in connection to holomorphic dynamics and multitangent functions ;  see also. \cite{Minh}, \cite{Minh4} and \cite{Murty-Sinha} for other articles dealing with these functions)~.

Such a function is an evaluation of a monomial quasi-symmetric function $M_{s_1, \cdots, s_r}$
(See at \cite{Gessel}, or \cite{Bergeron} for a more
recent presentation). Therefore, this family spans an algebra $\mathcal{H}mzf_{cv}$,
where the product is the product of $\textit{QSym}$, known as the stuffle product (also called
augmented shuffle, contracting shuffle, quasi-shuffle, $\cdots$ : see \cite{Hoffman}, \cite{Ebrahimi-FardGuo},
$\cdots$). The aim of this note is to describe the algebraic structure of $\mathcal{H}mzf_{cv}$~.

\bigskip

The familly of Hurwitz multizeta functions is a functional analogue of multizeta
values, which are defined by $\mathcal{Z}e^{s_1, \cdots, s_r} = \mathcal{H}e^{s_1, \cdots, s_r}(0)$
(See  \cite{Cartier}, \cite{Waldschmidt} and \cite{Zudilin}). Little is known on the 
arithmetical properties of  multizeta values, despite of the recent works \cite{Brown1}, \cite{Brown2} , \cite{Minh2} and \cite{Minh3}.
These numbers span an algebra denoted by $\mathcal{M}zv_{cv}$.
An algebraic description of it may be asked (for instance, if it is a graded algebra, and then a formula for
the dimensions of homogeneous components, etc). 

Elucidating the algebraic structure of an algebra, like $\mathcal{M}zv_{cv}$ or $\mathcal{H}mzf_{cv}$, is
in general a difficult question, since no one can be sure there does not exist a second natural product
hidden somewhere inside the algebraic structure itself. This is for instance the case in
$\mathcal{M}zv_{cv}$. In the case of Hurwitz multizeta functions, the situation is considerably
simpler, because the stuffle product, coming from the identification of the Hurwitz multizeta functions with a specialisation of monomial quasi-symmetric functions, is actually the only product to be taken in account.

In this note, we will explain the reasons which elucidate completely the case of $\mathcal{H}mzf_{cv}$.
The key point of the proof, which will be the heart of this note, is the following

\bigskip

\begin{Theorem}	\label{linear independance of HMZF}
	The family	$	\left (\mathcal{H}e^{s_1, \cdots, s_r} \right )
				_{	r \in \N
					\atop
					{	s_1, \cdots, s_r \in \N^*
						\atop
						s_1 \geq 2
					}
				}
			$ is $\C(z)$-linearly independent.
\end{Theorem}

Let us notice that \cite{Joyner}  proposes a proof of a weaker property, the linear independence over $\C$ of the Hurwitz multizeta functions\footnote{According to the referee who pointed out this reference, this proofs contains a few errors...}.

The following corollaries answer all the questions that one would ask for the multizeta values.

\bigskip

\begin{Corollary}
	Let us denote by $\textit{QSym}_{cv}$ the subalgebra of $\textit{QSym}$ spanned by the monomials
	$M_I$ with a composition $I = (i_1, \cdots, i_r)$ such that $i_1 \geq 2$~.
	Then:
	\begin{equation}
		\mathcal{H}mzf_{cv} \simeq \textit{QSym}_{cv}
				    \simeq
				    \Q	\big \langle
							\mathcal{L}yn(y_1 ; y_2 ; \cdots) - \{y_1\}
					\big \rangle
		\ ,
	\end{equation}
	where $\mathcal{L}yn(y_1 ; y_2 ; \cdots)$ denotes the set of Lyndon words over the alphabet
	$Y = \{ y_1 ; y_2 ; \cdots \}$~.
\end{Corollary}

For more information about Lyndon words, we refer the reader to \cite{Reutenauer}.

\bigskip

\begin{Corollary}
	Each algebraic relation in $\mathcal{H}mzv_{cv}$ comes from the expansion of
	stuffle products.
\end{Corollary}

\bigskip

\begin{Corollary}
	Let us denote the algebra of Hurwitz multizeta functions of weight $n$ by $\mathcal{H}mzv_{cv,n}$, 
	that is, the subalgebra of $\mathcal{H}mzv_{cv}$ spanned by the Hurwitz multizetas
	$\mathcal{H}e^{s_1, \cdots, s_r}$ such that $s_1 + \cdots + s_r = n$~.	\\
	Then:
	\begin{enumerate}
	      \item $\mathcal{H}mzv_{cv}$ is graded by the weight:
			\begin{equation}
				\left \{
					\begin{array}{l}
						\mathcal{H}mzv_{cv} = \bigoplus	_{n \in \N}
										\mathcal{H}mzv_{cv,n}
						\ .
						\\
						\forall (p , q) \in \N^2 \ , \ 
						\mathcal{H}mzv_{cv,p} \cdot \mathcal{H}mzv_{cv,q}
						\subset
						\mathcal{H}mzv_{cv,p + q} \ .
					\end{array}
				\right.
			\end{equation}
		\item $	\left \{
					\begin{array}{l}
						\text{dim } \mathcal{H}mzf_{cv,0} = 1
						\ , \ 
						\text{dim } \mathcal{H}mzf_{cv,0} = 0.
						\\
						\text{dim } \mathcal{H}mzf_{cv,n} = 2^{n - 1}
						\text{ for all }n \geq 2.
					\end{array}
				\right.
			$
 
	\end{enumerate}
\end{Corollary}

\section {A fundamental lemma}

Hurwitz multizeta functions are ``translations'' of multizeta values. Therefore, it is natural
to examine how the shift operator acts on such functions. It turns out that the Hurwitz multizeta
functions satisfy a fundamental difference equation.

\begin{Lemma}	
	Let $J^{s_1, \cdots, s_r}$, for $r \in \N$, be the function such that
	$	J^{s_1, \cdots, s_r} (z)
		=
		\left \{
			\begin{array}{cl}
				\displaystyle	{ \frac {1} {z^{s_1}} }	&	\text{, if } r = 1 \ .
				\vspace{0.1cm}	\\
				0					&	\text{, otherwise } \ .	\\
			\end{array}
		\right.
	$
	\\
	Then, for all $r \in \N$, $(s_1, \cdots, s_r) \in (\N^*)^r$ such that $s_1 \geq 2$, we have:
	\begin{equation}		\label{difference equation}
		\Delta_- (\mathcal{H}e^{s_1, \cdots, s_r})(z)
		=
		\frac {1} {z^{s_r}}
		\mathcal{H}e^{s_1, \cdots, s_{r - 1}}(z)
		\text { , where }
		\Delta_- (f) (z) = f (z - 1) - f(z) \ .
	\end{equation}
\end{Lemma}

\section {An example of linear independence of Hurwitz multizeta functions}

As an example, let us consider three rational functions $F_2$, $F_{2,1}$ and $G$
valued in $\C$ such that
\begin{equation}
	F_2 \mathcal{H}e^2 + F_{2,1} \mathcal{H}e^{2,1} = G \ .
\end{equation}
If $F_{2,1} \neq 0$, we can assume that $F_{2,1} = 1$. Then, applying the operator $\Delta_-$ to this
relation, we obtain:
\begin{equation} \label{simple equation}
	\widetilde{F}_2 \mathcal{H}e^2 = \widetilde{G} \ ,
\end{equation}
where
\begin{equation}
    \left \{
	\begin{array}{l}
	    \widetilde{F}_2 (z) = \Delta_- (F_2) (z) + \frac {1} {z} = \Delta_- (F_2 + \mathcal{H}e^1) (z)
	    \vspace{0.1cm}\\
	    \widetilde{G} (z) = \Delta_- (G) (z) - \frac {1} {z^2} \cdot F_2 (z - 1)
	\end{array}
    \right.
    \text{ and \ }
    \mathcal{H}e^1(z) = \sum	_{n > 0}
				\left (
					\frac	{1}	{n + z}
					-
					\frac	{1}	{n}
				\right )
    \ .
\end{equation}
From \eqref{simple equation}, we can easily deduce that $\widetilde{F}_2 = \widetilde{G} = 0$~.
Consequently, $F_2 + \mathcal{H}e^1$ is a $1$-periodic function. But, looking at the possible poles,
we can conclude that this is not possible. Therefore, $F_{2,1} = 0$, which also
implies $F_2 = G = 0$.

Thus, the function $1$, $\mathcal{H}e^2$ and $\mathcal{H}e^{2,1}$ are $\C(z)$-linearly independent.

\bigskip

This argument will be generalized in the proof of theorem \ref{linear independance of HMZF} and
can be enlightened by the following

\begin{Lemma}	\label{a contredire}
	Let $F$ be a rational fontion and $f$ a $1$-periodic function.			\\
	If, for a n-tuple $(\lambda_1 , \cdots , \lambda_n) \in \C^n$ ,
	$\displaystyle	{	F + \sum	_{i = 1}
									^n
									\lambda_i \mathcal{H}e^{i}
						=
						f
					}
	$ holds, then we necessarily have:
	$$	\left \{
				\begin{array}{l}
					\lambda_1 = \cdots = \lambda_n = 0\ .	\\
					F \text{ and } f \text{ are constant functions .}
				\end{array}
		\right.
	$$
\end{Lemma}

\section {Linear independence of Hurwitz multizeta functions over the rational function field}

We will give a sketch of proof of Theorem \ref{linear independance of HMZF} based on an induction
process. Let us introduce some notations, an order and properties:

\noindent
$a$. For all $d \in \N$, let $\mathcal{S}^{\star}_{\leq d}$ and $\mathcal{S}^\star_d$ be the sets defined
by
$\mathcal{S}^{\star}_{\leq d} = \left \{
					\seq{s} \in \mathcal{S}^\star \ ;\ d \text{\degre} \seq{s} \leq d
				\right \}
$ and	$\mathcal{S}^\star_d =	\left \{
					\seq{s} \in \mathcal{S}^\star \ ;\ d \text{\degre} \seq{s} = d
				\right \}
	$~, where	$	\mathcal{S}^\star =	\{	(s_1, \cdots, s_r) \in (\N^*)^r
								\ ; \ 
								s_1 \geq 2
							\}
			$ and	$	d \text{\degre} \seq{s} = d \text{\degre} (s_1, \cdots, s_r)
								= s_1 + \cdots + s_r - r
				$.

\noindent
$b$. We can order the sequences of $\mathcal{S}^\star_{d + 1}$ by numbering first the sequences
of length $1$, then those of length $2$, etc. For the following proof, let us notice that it will not be necessary
to precise the numbering within a length. So, we will consider the sets
$\mathcal{S}^\star_{d + 1} = \{ \seq{s}^n \ ; \ n \in \N^*\}$ and, for $n \in \N$,
$S_n = \{ \seq{s}^i \ ; \ 1 \leq i \leq n\}$~.

\noindent
$c$. We will finally consider the following properties:

$$	\begin{array}{ll}
		\mathcal{D}(d) :	&	\text{``the family }
						\left ( \mathcal{H}e^{\seq{s}} \right )
						_{\seq{s} \in \mathcal{S}^\star_{\leq d}}
						\text { is } \C(z) \text {-linearly independent.''}
		\\
		\mathcal{P}(d,n) :	&	\text{``the family }
						\left ( \mathcal{H}e^{\seq{s}} \right )
						_{\seq{s} \in \mathcal{S}^\star_{\leq d}}
						\bigcup
						\left ( \mathcal{H}e^{\seq{s}} \right )
						_{\seq{s} \in S_n}
						\text { is } \C(z) \text {-linearly independent.''}
	\end{array}
$$

In order to prove the theorem, we only need to give details on the heredity of the property
$\mathcal{D}(n)$~. This will also be done by an induction process and the proof boils down to
the following implication:
$$	\forall (d,n) \in \N^2 \ , \ \mathcal{P} (d,n) \Longrightarrow \mathcal{P}(d, n + 1)\ .$$

In order to do this, let us assume that Property $\mathcal{P}(d,n)$ holds for a pair
$(d,n) \in \N^2$ and consider the relation:
\begin{equation}	\label{dependance lineaire}
	\sum	_{	\seq{s} \in \mathcal{S}^\star_{\leq d} \cup S_{n + 1}
				\atop
				\seq{s} \neq \emptyset
			}
			F_{\seq{s}} \mathcal{H}e^{\seq{s}}
	=
	F
	\ ,
\end{equation}
where $F$ and $F_{\seq{s}}$ , $\seq{s} \in \mathcal{S}^\star_{\leq d} \cup S_{n + 1}$ , are
rational functions. Thanks to $\mathcal{P}(d,n)$, it is sufficient to prove that
$F_{\seq{s}^{n + 1}} = 0$~ . Using a proof by contradiction, we can then assume that
$F_{\seq{s}^{n + 1}} = 1$~.	\\

\noindent
\textbf{\underline{Step $1$ :}} Application of $\Delta_-$ to relation $(\ref{dependance lineaire})$~.
\\

Applying $\Delta_-$ to $(\ref{dependance lineaire})$ and
using the induction hypothesis $\mathcal{P}(d,n)$, we obtain 
\begin{equation}	\label{systeme}
	\left \{
		\begin{array}{l}
			\forall \seq{s} \in \left ( \mathcal{S}^\star_{\leq d} \cup S_n \right ) - \{\emptyset\}
			\ , \ 
			\displaystyle	{	\Delta_- (F_{\seq{s}})(z) + \sum	_{	k \in \N^*
												\atop
												\seq{s} \cdot k \in \mathcal{S}^\star_{\leq d} \cup S_n
											}
											F_{\seq{s} \cdot k}(z - 1) \cdot J^k(z)
						=
						0\ . 
					}
			\\
			\displaystyle	{	\Delta_-  (F)(z) = \sum	_{k = 2}
									^{d + 1}
									F_{k}(z - 1) J^k(z)
									+
									(1 - \delta_{n,0}) F_{d + 2}(z - 1) J^{d + 2}(z)\ .
							}
		\end{array}
	\right.
\end{equation}

\noindent
\textbf{\underline{Step $2$:}} A lemma which gives some partial solutions to the system $(\ref{systeme})$~.
\\

Denoting by $\seq{s}^{n + 1} = \seq{u} \cdot p$ with $p \geq 1$ and
$\seq{u} \in \mathcal{S}^\star_{\leq d} \cup S_n$, we can prove that:

\begin{Lemma}	\label{partial system's resolution}
	Let $r$ be a positive integer and $p \geq 2$~.	\\
	Let us also consider two $r$-tuples, $(n_1 ; \cdots ; n_r) \in \N^r$ and $(k_1 ; \cdots ; k_r) \in (\N^*)^r$ such that
	$	\displaystyle	{	\sum	_{i = 1}
						^r
						(k_i - 1)
					\leq
					p - 2
				}
	$~.
	Then,	$F_{\seq{u} \cdot k_1 \cdot 1^{[n_1]} \cdots \cdot k_r \cdot 1^{[n_r]}}$ is null
	if $n_r > 0$, and is constant if $n_r = 0$~.
\end{Lemma}

\bigskip

\noindent
\textbf{\underline{Step $3$:}} Revealing the contradiction.
\\

The system $(\ref{systeme})$, applied to $\seq{u}$, gives us
\begin{equation}
	\Delta_-	\left (
				F_{\seq{u}}	+ \sum	_{k = 1}
							^{p - 1}
							F_{\seq{u} \cdot k} \mathcal{H}e^k
						+ \mathcal{H}e^p
			\right )
	=
	0
	\ ,
\end{equation}
where $F_{\seq{u} \cdot 1}$ , $\cdots$ , $F_{\seq{u} \cdot p - 1}$ are rational fractions
which turn out to be constant according to the lemma \ref{partial system's resolution}
and now denoted by $f_{\seq{u} \cdot 1}$ , $\cdots$ , $f_{\seq{u} \cdot p - 1}$. Thus,
$	\displaystyle	{	F_{\seq{u}}
				+
				\sum	_{k = 1}
					^{p - 1}
					f_{\seq{u} \cdot k} \mathcal{H}e^k
				+
				\mathcal{H}e^p
			}
$ defines a $1$-periodic function. The coefficients being not all zero in this relation,
this contradicts Lemma \ref{a contredire} and concludes the proof.

\end{document}